\documentclass{amsart}
\usepackage{amsmath,amsthm,amssymb,amsfonts}
\usepackage[colorlinks=true]{hyperref}

\hypersetup{urlcolor=blue, citecolor=blue}

\def\doi#1{   {\href{http://dx.doi.org/#1}
   {{\mdseries\ttfamily DOI}}}}

\newcommand{\al}{\alpha}    \newcommand{\be}{\beta}
    
  \newcommand{\ep}{\varepsilon}

\newcommand{\R}{\mathbb{R}}

\newcommand{\pt}{\partial_t}\newcommand{\pa}{\partial}
\newcommand{\les}{{\lesssim}}
\newcommand{\beeq}{\begin{equation}}\newcommand{\eneq}{\end{equation}}

\newcommand{\cd}{{\,\cdot\,}}
\newcommand{\Sp}{{\mathbb S}}\def\CO{\mathcal {O}}

\newenvironment{prf}{\noindent {\bf Proof.} }{\endprf\par}
\def \endprf{\hfill  {\vrule height6pt width6pt depth0pt}\medskip}

\numberwithin{equation}{section}

\newcommand{\gm}{\mathfrak{g}}

\def\<{\langle}             \def\>{\rangle}
\def\({\left(}                 \def\){\right)}

\newtheorem{theorem}{Theorem}[section]

\newtheorem{lemma}[theorem]{Lemma}

\theoremstyle{definition}

\theoremstyle{remark}\newtheorem{rem}{Remark}

\theoremstyle{definition}

\title[Global existence of semilinear damped wave equations]
      {
      Global Existence for semilinear  damped wave equations in relation with the Strauss conjecture
}

\author{Mengyun Liu}
\address{School of Mathematical Sciences\\                Zhejiang University\\                Hangzhou 310027, P. R. China}
\email{mengyunliu@zju.edu.cn}

\author{Chengbo Wang}\address{School of Mathematical Sciences\\                Zhejiang University\\                Hangzhou 310027, P. R. China}\email{wangcbo@zju.edu.cn}
\urladdr{http://www.math.zju.edu.cn/wang}

\date{\today}
\dedicatory{} \commby{}

\begin{document}
\bibliographystyle{plain}

\begin{abstract}
We study the global existence of solutions to semilinear wave equations with power-type nonlinearity and general lower order terms on $n$ dimensional nontrapping asymptotically Euclidean manifolds, when $n=3, 4$. In addition, we prove almost global existence with sharp lower bound of the lifespan for the four dimensional critical problem.
\end{abstract}



\maketitle

\section{Introduction}
Let $(\R^{n}, \gm)$ be a nontrapping asymptotically Euclidean (Riemannian) manifolds, with $n\ge 3$,
\begin{align}
\label{ea6}
&\gm=g^{0}+g^1(r)+g^2(x), \gm \mathrm{\ is  \ nontrapping},
\end{align}
where $\gm=\gm_{ij}dx^i dx^j$ and $\gm_{ij}=g^0_{ij}+
g^1_{ij}+ g^2_{ij}$,
$g^0_{ij}= \delta_{ij}$,
$(\gm^{ij}(x))$ denotes the inverse matrix of $(\gm_{ij}(x))$.
Here and in what follows, the Einstein summation convention is performed over repeated upper and lower indices,  $1\le i,j\le n$.
 We assume the first perturbation $g^1$ is radial and for some fixed $\rho_{1}>0, \rho_{2}>1$,
\begin{align}
\label{ea8}
&|\nabla_x^a g^l_{ij}|\les_a\langle x\rangle^{-|a|-\rho_l}, l=1,2, \rho_1<\rho_2,\rho=\min(\rho_{1}, \rho_{2}-1) > 0, |a|\le 4,
\end{align}
where $\langle x\rangle = \sqrt{1 + x^{2}}$. On $(\R^n, \gm)$ with $n\ge 3$, it is known from the works
of Bony-H\"afner \cite{BoHa} and Sogge-Wang \cite{SW10} that we have the local energy estimates without loss and Strichartz estimates (see also Metcalfe-Sterbenz-Tataru \cite{MST}).
The nonlinear wave equations on the nontrapping asymptotically Euclidean manifolds has received much attention in recent years. For example, Bony-H\"afner \cite{BoHa}, Sogge-Wang \cite{SW10} and Yang \cite{Yang13} studied the analogs of the John-Klainerman theorem \cite{JK84} and global existence under null conditions for semilinear wave equations (see also Wang-Yu \cite{WY14} and Yang \cite{Yang16} for quasilinear wave equations).
Sogge-Wang \cite{SW10}, Wang-Yu \cite{WaYu11},
Metcalfe-Wang \cite{MW17} and Wang \cite{W17}
proved the analogs of the global existence part of the Strauss conjecture when $n=3,4$ (see also \cite{Wang18} for a review of recent results and
Wakasa-Yordanov
\cite{WaYo18-1}
for the recent blow-up results with critical power). For the analogs of the Glassey conjecture, see Wang \cite{W15} and references therein.

In this paper, we are interested in the small data global existence of solutions for the Cauchy problem of the following semilinear wave equations with general lower order term, in relation with the Strauss conjecture, posed on nontrapping asymptotically Euclidean (Riemannian) manifolds
\begin{equation}
\label{1.1}
\begin{cases}
u_{tt} - \Delta_{\gm}u + \mu(t,x)\pa_{t}u+
\mu^{j}(t,x)\pa_{j}u + \mu_0(t,x)u  = F_{p}(u),\\
u(0,x)=u_0(x),  \pt{u}(0,x)=u_1(x).
\end{cases}
\end{equation}
Here,
$\Delta_{\gm} =
\sqrt{|\gm|}^{-1}\partial_{i}g^{ij}\sqrt{|\gm|}\partial_{j}$ is the standard Laplace-Beltrami operator, with
$|\gm|=\det(\gm_{ij}(x))$. The nonlinearity $F_p$ is assumed to behave like $|u|^p$, more precisely, we assume $$
|u|\ll 1\Rightarrow
\sum_{0\leq l\leq2}|u|^{l}|\pa_{u}^{l}F_{p}(u)| \les |u|^{p}\ .$$
Concerning lower order terms, we assume
 \beeq
 \label{ea66}
 Y^{\leq 2}(\mu, \mu^{j}) \in L^{1}_{t}(L_{x}^{\infty}\cap \dot{W}^{1,n}),
  |x| Y^{\leq 2}\mu_0 \in L^{1}_{t}(L_{x}^{\infty}\cap \dot{W}^{1,n}).
 \eneq

When $\gm = g^{0}$ and $\mu=\mu^{j} = \mu_0 = 0$, the problem has been extensively investigated and is known as the Strauss conjecture,
which was initiated in the work of John \cite{John79}.
 It is known that, in general (with $F_p(u)=|u|^p$),
the problem admits small data global existence only if $p> p_c(n)$ (see Yordanov-Zhang \cite{YorZh06} and
Zhou \cite{Zh07}), where
the critical power $p_{c}(n)$ is the positive root of equation \beeq\label{eq-Strauss}(n-1)p^{2} - (n+1)p - 2 = 0\ .\eneq
The global existence for
 small initial data when $p\in (p_c(n), 1+4/(n-1))$  followed in Georgiev-Lindblad-Sogge
\cite{GLS97} (see also Tataru
\cite{Ta01-2}). See
Wang-Yu \cite{WaYu12rev}, Wang \cite{Wang18} for a complete history and recent works for the problem on various space-time manifolds.

On the other hand, there are many recent works concerning the damped wave equations, $\gm = g^{0}$ and $\mu^{j} = \mu_0 = 0$, with typical damping term depending only on time $$\mu(t,x)=\frac{\mu}{(1+t)^{\beta}}, \mu>0\ .$$
For the case $\be<1$, the damping term is strong enough to make the  problem behaves like heat equations
and the problem has been well-understood.
There are some interesting critical phenomena happening for the scale-invariant case $\be=1$ and it appears that the critical power is $p_c(n+\mu)$ for relatively small $\mu>0$.

For the remaining case, $\be>1$ (which is also referred as the scattering case), where the damping term is integrable, it is natural to expect that the problem behaves like the nonlinear wave equations without damping term. In a recent work of Lai-Takamura \cite{LaiTa18}, the authors proved blow up results for $1<p<p_c(n)$, together with upper bound of the lifespan.
In particular, it is shown that for $1<p<p_c(n)$ with $n\ge 2$, we have
\beeq\label{eq-life-subcrit}
T_\ep\le C \ep^{\frac{2p(p-1)}{(n-1)p^{2} - (n+1)p - 2}},
\eneq
where $T_\ep$ denotes the lifespan
and $\ep$ is the size of the initial data.
For the critical case, $p=p_c(n)$,
with $\gm=g^{0}+g^{2}(x)$ and general (nonnegative) damping term $\mu=\mu(t)\in L^1$, where
$|g^{2}_{ij}|+|\nabla g^{2}_{ij}|\le C e^{-\al (1+|x|)}$ for some $\al>0$, at the final stage of preparation of the current manuscript, we learned that Wakasa-Yordanov \cite{WaYo18-2} obtained the expected exponential upper bound of the lifespan,
\beeq\label{eq-life-crit}
T_\ep\le \exp(C\ep^{-p(p-1)})\ .
\eneq

In this paper, we are interested in
complementing to the blow up results for the scattering case, by proving global existence results
on general nontrapping asymptotically Euclidean  manifolds. Moreover, in the process, we find that we could handle more general lower order perturbation terms as in \eqref{1.1}, under the assumption \eqref{ea66}. The first main theorem of this paper states as follows:
\begin{theorem}
\label{main.0}
Let $n =3, 4$, and
consider the Cauchy problem \eqref{1.1} with
\eqref{ea66}, posed on
nontrapping asymptotically Euclidean manifolds, with
\eqref{ea6} and \eqref{ea8}. Then if $p > p_{c}(n)$,
the problem \eqref{1.1} admits a unique global solution for any initial data which are sufficiently small, decaying and regular.
\end{theorem}

For more precise statement, see Theorem \ref{global}.
Concerning the proof, the idea is to adapt the recent approach of using local energy and weighted Strichartz estimates, which has been very successful in the recent resolution of the Strauss conjecture on various space-times, including Schwarzschild/Kerr black-hole space-times (\cite{FW11}, \cite{HMSSZ}, \cite{LMSTW}, \cite{MW17}). In particular, we revisit the proof of \cite[Theorem 4.1]{MW17} to extract the key weighted Strichartz estimates, Lemma \ref{3.3}, which, combined with the local energy estimates (\cite{BoHa},  \cite{SW10}, \cite{W15}, \cite{MST}), is good enough to treat the lower order terms in  \eqref{1.1} in a perturbative way.

When there are no global solutions, it is also interesting to obtain sharp estimates of the lifespan.
On this respect,
it turns out that our argument could also be adapted to show some of the sharp results, as long as we could add some favorable terms in the desired space-time estimates, which could be used to absorb the lower order term in  \eqref{1.1}.
To illustrate the argument, as an example, we prove the following lower bound estimate of the lifespan for the four-dimensional critical problem, which is sharp in general, comparing with \eqref{eq-life-crit} of \cite{WaYo18-2}.
\begin{theorem}
\label{main.1}
Let $n =4$, $F_p=u^2$ and
consider the Cauchy problem \eqref{1.1} with
\beeq
\label{6666}
 Y^{\leq 3}(\mu, \mu^{j}),
  |x| Y^{\leq 3}\mu_0 \in L^{1}_{t}(L_{x}^{\infty}\cap \dot{W}^{1,4}),
  \eneq
   posed on
nontrapping asymptotically Euclidean manifolds, with
\eqref{ea6} and \eqref{ea8}. Then there exists $c>0$,
such that the problem \eqref{1.1} admits almost global solution, up to
\beeq\label{eq-life-crit-2}
T_\ep= \exp(c\ep^{-2})\ ,
\eneq
 for any initial data which are sufficiently small (of size $\ep\ll 1$), decaying and regular.
\end{theorem}
See Theorem \ref{almost} for more precise statement.

\begin{rem}
It is remarkable that, in our statement, the damping coefficient $\mu(t,x)$ is not required to be nonnegative, which were assumed in
both \cite{LaiTa18} and
\cite{WaYo18-2}. Moreover, the authors believe that the  nonnegative assumption there are not necessary.
\end{rem}
\begin{rem}
On the two dimensional Euclidean space, i.e., $n=2$ with $g^{1} = g^{2} =0$, the idea in this paper can be exploited to show that the weighted Strichartz estimates of
\cite{FW11}
and \cite{HMSSZ} are strong enough to yield
small data global existence for \eqref{1.1} with $p > p_{c}(2)$.
\end{rem}
\begin{rem}
Based on
the existence results in previous works, as we have illustrated in our theorems, our argument could be adapted to show the following
 lower bounds of the lifespan  for \eqref{1.1}
 \begin{enumerate}
  \item $n=2$, $2<p<p_c(2)$, $\gm=g^0$, $T_\ep\ge c  \ep^{\frac{2p(p-1)}{ p^{2} - 3 p - 2}}$, \cite[Theorem 6.1]{JWY12}.
  \item $n=2$, $p=p_c(2)$, $\gm=g^0$, $T_\ep\ge \exp(c  \ep^{-(p-1)^2/2})$, \cite[Theorem 6.2]{JWY12}.
  \item $n=3$, $2\le p<p_c(3)$,  $T_\ep\ge c  \ep^{\frac{p(p-1)}{ p^{2} - 2 p - 1}}$, \cite[Theorem 3.2]{W17}.
    \item $n=3$, $p=p_c(3)$,
    $T_\ep\ge \exp(c  \ep^{-2(p-1)})$, \cite[Theorem 3.2]{W17}.
\end{enumerate}
These lower bounds, together with the upper bounds available from \cite{LaiTa18}, show the sharpness of the lifespan estimates for $n=2$ with $p\in (2, p_c(2))$,
$n=3$ with $p\in [2, p_c(3))$.
\end{rem}



\subsection{Notation}
The vector fields to be used will be labeled as
\[Y=(Y_1,\dots, Y_{n(n+1)/2}) = \{\nabla_x,\Omega\}\ ,\]
where $\Omega$ denotes the
generators of spatial rotations
$\Omega_{ij} = x_i\partial_j - x_j\partial_i$, $1\le i<j\le n$.
For a norm $X$ and a nonnegative integer $m$, we shall use the shorthand
\[|Y^{\le m} u| = \sum_{|\mu|\le m} |Y^\mu u|,\quad \|Y^{\le m} u\|_X
= \sum_{|\mu|\le m} \|Y^\mu u\|_X.\]
Let $L^q_\omega$ be the standard Lebesgue space on the sphere
$\Sp^{n-1}$, we will use the following convention for mixed norms
$L^{q_1}_tL^{q_2}_rL^{q_3}_\omega$, with $r=|x|$ and $\omega\in \Sp^{n-1}$:
\[\|f\|_{L^{q_1}_tL^{q_2}_rL^{q_3}_\omega} = \Bigl\|\Bigl(\int
\|f(t,r\omega)\|^{q_2}_{L^{q_3}_\omega}
r^{n-1}\,dr\Bigr)^{1/q_2}\Bigr\|_{L^{q_1}(\{t\ge 0\})},\]
with trivial modification for the case $q_2=\infty$. Clearly $L_{t}^{q_{1}}L_{x}^{q_{2}} = L^{q_1}_tL^{q_2}_rL_{\omega}^{q_2}$. We denote $L^{q_{1}}_{T}L_{x}^{q_{2}} = L^{q_{1}}([0, T]; L_{x}^{q_{2}}(\R^n))$ for some $T >0$. Occasionally,
when the meaning is clear, we shall omit the subscripts.
As usual, we use $\|\cdot\|_{E_m}$ to denote the energy norm of order $m\ge0$,
\begin{align}
\label{1.8}
&\|u\|_E=\|u\|_{E_0}=\|\partial u\|_{L_t^\infty L_x^2},\ \|u\|_{E_m}=\sum_{|a|\leq m}\|Y^au\|_E.
\end{align}
We will use $\|\cdot\|_{LE}$ to denote the (strong) local energy norm
\begin{align}
\label{1.9}
&\|u\|_{LE}=\|u\|_E+\|\partial u\|_{\ell_\infty^{-1/2}(L_t^2L_x^2)}+\|r^{-1}u\|_{\ell_\infty^{-1/2}(L_t^2L_x^2)}\,
\end{align}
and $\|u\|_{LE_m}=\sum_{|a|\leq m}\|Y^au\|_{LE}$, where we write
$$\|u\|_{\ell_{q}^{s}(A)} = \|\phi_{j}(x)u(t,x)\|_{\ell_{q}^{s}(A)}
= \|\left(2^{js}\|\phi_{j}(x)u(t,x)\|_{A}\right)\|_{\ell^{q}(j\ge 0)},
$$
for a partition of unity subordinate to the (inhomogeneous) dyadic annuli, $\Sigma_{j \geq 0}\phi_{j}^{2}=1$. We denote $\|u\|_{LE_{T}}=\|u\|_{E_{T}}+\|\partial u\|_{\ell_\infty^{-1/2}(L_T^2L_x^2)}+\|r^{-1}u\|_{\ell_\infty^{-1/2}(L_T^2L_x^2)}$ and $\|u\|_{E_{T}} =\|\partial u\|_{L_T^\infty L_x^2}$.

\section{Preliminary}
In this section, we collect some inequalities we shall use later.
\begin{lemma}
Let $s \in [0,1]$ and $n\ge 3$. If $f\in L^{\infty}(\R^{n}) \cap \dot{W}^{1,n}(\R^{n})$, $g \in \dot{H}^{s-1}$, then we have 
\beeq
\label{2.1}
\|fg\|_{\dot{H}^{s-1}} \les \|f\|_{L^{\infty}\cap \dot{W}^{1,n}}\|g\|_{\dot{H}^{s-1}}.
\eneq
\end{lemma}
\begin{prf}
When $s=1$, by H\"{o}lder's inequality
$$\|fg\|_{L^{2}} \les \|f\|_{L^{\infty}}\|g\|_{L^{2}}.$$
When $s=0$, by duality, we need only to show
$$\|fg\|_{\dot{H}^{1}} \les \|f\|_{L^{\infty}\cap \dot{W}^{1,n}}\|g\|_{\dot{H}^{1}}.$$
Since $\|fg\|_{\dot{H}^{1}}\le \|\nabla f g\|_{L^{2}} + \|f \nabla g\|_{L^{2}}$, by  H\"{o}lder's inequality and Sobolev embedding
\begin{align*}
\|\nabla f g\|_{L^{2}} + \|f \nabla g\|_{L^{2}}&\les \|\nabla f\|_{L^{n}}\|g\|_{L^{\frac{2n}{n-2}}}+ \|f\|_{L^{\infty}}\|\nabla g\|_{L^{2}}\\
&\les \|f\|_{L^{\infty}\cap \dot{W}^{1,n}}\|g\|_{\dot{H}^{1}}.
\end{align*}
By interpolation, (\ref{2.1}) follows.
\end{prf}
\begin{lemma}
Let $s \in [0,1]$ and $n\ge 3$, then we have
\beeq
\label{2.2}
\|u/r\|_{\dot{H}^{s-1}} \les \|u\|_{\dot{H}^{s}}.
\eneq
\end{lemma}
\begin{prf}
When $s=1$, it is the classical Hardy's inequality
$$\|u/r\|_{L^{2}} \les \|u\|_{\dot{H}^{1}}\ .$$
By duality, we get
$$\|u/r\|_{\dot{H}^{-1}} \les \|u\|_{L^{2}},$$
which is the case $s=0$. Then by interpolation, (\ref{2.2}) follows.
\end{prf}
\section{Space-time estimates}
In this section, we collect various space-time estimates for linear wave equation
\beeq
\label{3.1}
u_{tt} - \Delta_{\gm}u = F(=F_{1} + F_{2})\ .
\eneq

\begin{lemma}[Local energy estimates]
\label{3.2}
Let $n\ge 3$ and consider linear wave equation \eqref{3.1} satisfying \eqref{ea6}, \eqref{ea8}.
Then
we have the following higher order local energy estimates
\begin{eqnarray}
\|u\|_{LE_m}\les \|\partial Y^{\le m}u(0)\|_{L_x^2}+\|Y^{\le m}F\|_{L_t^1L_x^2}, 0\le m\le 3\ .
\end{eqnarray}
\end{lemma}
\begin{prf}
It is proven in Wang \cite[Lemma 3.5]{W15}. The result with $g^1=0$ has been proven in Bony-H\"afner \cite{BoHa} and Sogge-Wang \cite{SW10}. See also Metcalfe-Sterbenz-Tataru \cite{MST}.
\end{prf}
\begin{lemma}[Weighted strichartz estimates]
\label{3.3}
Let $n\ge 3$. Consider linear wave equation \eqref{3.1} with \eqref{ea6}, \eqref{ea8}.  Then there exists $R>0$ so that if
$\psi_R$ is identically $1$ on $B_{2R}^c$ and vanishes on $B_R$, for $0\leq m \leq 2$,
we have
\begin{align}
  \label{2.3}
&  \|\psi_R r^{\frac{n}{2}-\frac{n+1}{p}-s} Y^{\le m}u\|_{L^p_{t,r}
    H_\omega^{\sigma}}+ \|\pa(\psi_RY^{\le m}u)\|_{L_{t}^{\infty}\dot{H}^{s-1}} \\\lesssim &
  \|\psi_R  Y^{\le m} u(0,\cd)\|_{\dot{H}^s} + \| \psi_R Y^{\le m} \partial_t  u(0,\cd)\|_{\dot{H}^{s-1}}\nonumber
\\&+
  \|\psi_R^p r^{-\frac{n-2}{2}-s}Y^{\le m}F_{1} \|_{L^1_{t,r}H^{s-1/2}_\omega} +
  \|Y^{\le m} F_{1}\|_{L^1_t L^2_x}+\|\psi_{R}Y^{\le m}F_{2}\|_{L^{1}_{t}\dot{H}^{s-1}}\nonumber
\end{align}
  for any $p\in (2,\infty)$, $s\in (1/2-1/p,1/2)$, and $0\le \sigma< \min(s-1/2+1/p,1/2-1/p)$.
\end{lemma}
\begin{prf}
It is essentially proved in \cite[Theorem 4.1]{MW17}.
\end{prf}

\begin{lemma}[Space-time estimates]\label{thm-wStri-homo}
Let $n\ge 4$. Consider linear wave equations \eqref{3.1} with \eqref{ea6}, \eqref{ea8}.
Then there exist a $R \gg 1$ so that $\psi_R$ is identically $1$ on $B_{2R}^c$ and vanishes on $B_R$, for any $T>2$ and $m \in [0, 3]$, we have
\begin{align*}
&~~~~~~~(\ln T)^{-1/4}\|\psi_{R}\<r\>^{-\frac{n-3}{4}} Y^{\leq m}u\|_{L_T^{4}L^{2}}+\|Y^{\leq m}u\|_{LE_{T}}\\
&~~~~~~~+\|\langle r\rangle^{-\frac{n-1}{4}}Y^{\leq m}u\|_{L^{2}_T L^{2}}
+\|\pa (\psi_{R}Y^{\leq m}u)\|_{L_{T}^{\infty}\dot{H}^{-1}}\\
&\les \|\pa^{\le 1} Y^{\leq m}u(0,\cdot)\|_{L^{2}}+\|\psi_{R}Y^{\leq m}\pa_{t}u(0)\|_{\dot{H}^{-1}}
+\|\psi_{R}^{2} r^{-(n-3)/2}Y^{\leq m}F_{1}\|_{L_{T}^2L^1H^{-1/2+}_{\omega}}\\
&+\|Y^{\leq m}(F_{1},F_{2})\|_{L^{1}_{T}L^{2}}+\|\psi_{R}Y^{\leq m}F_{2}\|_{L_{T}^{1}\dot{H}^{-1}}
\end{align*}
\end{lemma}
\begin{prf}
This is essentially proved in \cite[Lemma 5.6]{W17}, which is based on the local energy estimates
on nontrapping asymptotically Euclidean (Riemannian) manifolds, Lemma \ref{3.2}, as well as a sharper version of local energy estimates for small metric perturbation, due to Metcalfe-Tataru \cite{MeTa12MA}.
\end{prf}

\section{Global Existence}
In this section, we prove the global existence results,
Theorem \ref{main.0}.

\begin{theorem}
\label{global}
Let $n =3, 4$, and assume that \eqref{ea6}, \eqref{ea8} hold. Consider \eqref{1.1} with \eqref{ea66} and $p > p_{c}$. Set $s= \frac{n}{2}-\frac{2}{q-1} \in (\frac{1}{2}-\frac{1}{q}, \frac{1}{2})$ with $q =p$ if $p \in (p_{c}, 1+4/(n-1))$ and $q\in (p_{c}, 1+4/(n-1))$ is any fixed choice when $p \ge 1+4/(n-1)$. Then there exist $\ep_{0}$ sufficiently small and a $R$ sufficiently large, so that if $0 < \ep < \ep_{0}$ and
\beeq
\label{eqn-4.1}
\|\nabla^{\le 1}Y^{\le 2}u_{0}\|_{L^{2}_{x}}+ \|Y^{\le 2}u_{1}\|_{L^{2}_{x}}+\|\psi_{R}Y^{\le 2}u_{1}\|_{\dot{H}^{s-1}} \leq \ep,
\eneq
then there exists a global solution $u \in C([0, \infty);H^{3})\cap C^{1}([0, \infty);H^{2})$.
\end{theorem}
\begin{prf}
Without loss of generality, we  may assume $p \in (p_{c}, 1+4/(n-1))$. If not, one need only fix any $q \in (p_{c}, 1+4/(n-1))$ and apply the proof below while noting that Sobolev embeddings provide $\|u\|_{L^{\infty}_{t,x}} \les
\|u\|_{X_{2}}$, which suffices to handle the $p-q$ extra copies of the solution in the nonlinearity.\par
For $0 \leq m \leq 2$, we shall apply Lemma \ref{3.3} with $s= \frac{n}{2} - \frac{2}{p-1}$ and note that $s \in (\frac{1}{2}-\frac{1}{p}, \frac{1}{2})$ precisely $p_{c} < p < 1+4/(n-1)$. We set $-\al=\frac{n}{2}- \frac{n+1}{p} -s = \frac{2}{p-1} - \frac{n+1}{p}$, then $-\frac{n-1}{2}-s = -\al p$.
Let $\theta$ be a fixed number satisfying
$$2 < \theta < \min\(p, \frac{2(n-1)}{n-1-2\min\(s-\frac{1}{2}+\frac{1}{p}, \frac{1}{2}-\frac{1}{p}\)}\),$$
we define the norms
\begin{align*}
\|u\|_{X_{m}} &= \|r^{-\al}\psi_{R}Y^{\le m}u\|_{L^{p}L^{p}L^{\theta}}+\|\pa(\psi_{R}Y^{\le m}u)\|_{L_{t}^{\infty}\dot{H}^{s-1}}\\
&~~~+\|Y^{\le m}u\|_{\ell_{\infty}^{-\frac{3}{2}}L^{2}L^{2}L^{2}}+\|\pa Y^{\le m} u\|_{L^{\infty}L^{2}L^{2}}\ ,
\end{align*}
$$\|F\|_{N_{m}}=\|r^{-\al p}\psi_{R}^{p}Y^{\le m}F_{1}\|_{L^{1}L^{1}L^{2}}+ \|Y^{\le m}(F_{1},F_{2})\|_{L_{t}^{1}L_x^{2}}+\|\psi_{R}Y^{\le m}F_{2}\|_{L_{t}^{1}\dot{H}^{s-1}}$$
where $F = F_{1} + F_{2}$ and $R$ is the large constant occurred in Lemma \ref{3.3}.
Then by Lemma \ref{3.2} and Lemma \ref{3.3}, for linear wave equation (\ref{3.1}), there exist  constants $C_{0}, C_1>0$ such that, for any $0\le m\le 2$,
\beeq
\label{3.4}
\|u\|_{X_{m}} \leq C_{0}\big(\|\pa^{\leq 1} Y^{\le m}u(0,\cd)\|_{L^2}  +\|\pa \psi_R  Y^{\le m} u(0,\cd)\|_{\dot{H}^{s-1}}  + \|F\|_{N_{m}}\big)
\le C_1\ep+C_0\|F\|_{N_m}
\ ,\eneq
where we have used \eqref{eqn-4.1} for the initial data.
We set $u^{(0)} =0$ and recursively define $u^{(k+1)}$ be the solution to the linear equation
\begin{equation}
\label{4.1}
\begin{cases}
\pa_{tt}u^{(k+1)}-\Delta_{\gm}u^{(k+1)} +
 \mu \pt u^{(k+1)} +
 \mu^{j}\pa_{j}u^{(k+1)} + \mu_0u^{(k+1)}  = F_{p}(u^{(k)}),\\
u^{(k+1)}(0,x)=u_{0},  \pt{u^{(k+1)}}(0,x)=u_{1}.
\end{cases}
\end{equation}
We will prove that $u^{(k)}$ is well defined, bounded in $X_{2}$ and convergent in $X_{0}$.

\emph{\textbf{Well defined}}: It is easy to see $u^{(1)} \in C([0,\infty);H^{3})\cap C^{1}([0,\infty);H^{2})$. If we have $u^{(k)} \in C([0,\infty);H^{3})\cap C^{1}([0, \infty);H^{2})$ for some $k\geq 1$.
Since $$|\nabla^{\leq 2}F_{p}(u^{(k)})| \lesssim |u^{(k)}|^{p-1}|\nabla^{\leq 2}u^{(k)}| + |u^{(k)}|^{p-2}(\nabla^{\leq 1}u^{(k)})^{2},$$
by H\"{o}lder's inequality and Sobolev embedding, we have, for any $t\in [0,T]$ with fixed $T<\infty$,
\begin{align*}
\|\nabla^{\leq 2}F_{p}(u^{(k)})\|_{L^{2}} \les \|u^{(k)}\|^{p-1}_{L_{x}^{\infty}}\|\nabla^{\leq 2}u^{(k)}\|_{L^{2}}
+ \|u^{(k)}\|^{p-2}_{L_{x}^{\infty}}\|\nabla^{\leq 1}u^{(k)}\|^{2}_{L^{4}} \les \|u^{(k)}\|^{p}_{H^{3}} <\infty.
\end{align*}
Thus $F_{p}(u^{(k)}) \in L^{1}_{t}([0,T];H^{2})$ for any $T < \infty$. By standard existence theorem of linear wave equations (see, e.g., Sogge \cite{So08}), $u^{(k+1)} \in C([0,T];H^{3})\cap C^{1}([0,T];H^{2})$  for any $T < \infty$ and so
$u^{(k+1)} \in C([0,\infty);H^{3})\cap C^{1}([0, \infty);H^{2})$. Hence the iteration sequence is well defined.\par
\emph{\textbf{Boundedness}}: For $k =0 $, we rewrite the equation (\ref{4.1}) as
\begin{equation}
\label{4.7}
\begin{cases}
\pa_{tt}u^{(1)} - \Delta_{\gm}u^{(1)}  =F = -\mu \pt u^{(1)} - \mu^{j}\pa_{j}u^{(1)} - \mu_0 u^{(1)} \\
u^{(1)}(0,x)=u_{0},  \pt u^{(1)}(0,x)=u_{1}.
\end{cases}
\end{equation}
By applying (\ref{3.4}) to (\ref{4.7}) with $F_{1} = 0, F_{2} = -\mu \pt u^{(1)} - \mu^{j}\pa_{j}u^{(1)} - \mu_0 u^{(1)}$,  we have
\beeq
\label{4.2}
\|u^{(1)}\|_{X_{2}}\leq C_{1}\ep + C_{0}(\|\psi_{R} Y^{\le 2}F_{2}\|_{L^{1}_{t}\dot{H}^{s-1}}+\|Y^{\leq 2}F_{2}\|_{L_t^{1}L_x^{2}}).\eneq
By applying (\ref{2.1}) and (\ref{2.2}) to the last two terms, we have
\begin{align}
\|\psi_{R}&Y^{\le 2}F_{2}\|_{\dot{H}^{s-1}}+\|Y^{\leq 2}F_{2}\|_{L_x^{2}}
\les\|(Y^{\leq 2}(\mu,\mu^{j}), rY^{\leq 2}\mu_{0})\|_{L_{x}^{\infty}\cap \dot{W}^{1,n}}
\label{eq-4.1}\\
 &\times(\|\psi_{R} \pa Y^{\leq 2} u^{(1)}\|_{\dot{H}^{s-1}}+\|\pa Y^{\leq 2}u^{(1)}\|_{L_x^2}+\|\pa(\psi_{R} Y^{\leq 2}u^{(1)})\|_{\dot{H}^{s-1}}).\nonumber
\end{align}
Notice that
$\nabla_x \psi_{R}$ is compactly supported and $n \geq 3$, and so
\begin{align}
\nonumber
\|\psi_{R} \pa Y^{\leq 2} u\|_{\dot{H}^{s-1}} &\les \|\pa(\psi_{R} Y^{\leq 2}u)\|_{\dot{H}^{s-1}}+\|(\nabla_x \psi_{R})Y^{\leq 2}u\|_{\dot{H}^{s-1}}\\ \label{liu7}
&\les\|\pa(\psi_{R}Y^{\leq2}u)\|_{\dot{H}^{s-1}}+\|(\nabla_x \psi_{R})Y^{\leq2}u\|_{L^{\frac{2n}{n-2}}}\\ \nonumber
&\les \|\pa(\psi_{R}Y^{\leq2}u)\|_{\dot{H}^{s-1}}+\|\nabla Y^{\leq2}u\|_{L_x^{2}},
\end{align}
where we have used Sobolev embedding $\dot{H}^{1} \hookrightarrow L^{\frac{2n}{n-2}}$. Hence we obtain
from \eqref{4.2} and \eqref{eq-4.1} that
\begin{align}
 \|u^{(1)}\|_{X_{2}}&
\leq  C_{1}\ep
+C_{2}\int_{0}^{T}\|(Y^{\leq 2}(\mu,\mu^{j}), rY^{\leq 2}\mu_{0})\|_{L_{x}^{\infty}\cap \dot{W}^{1,n}}\label{eqn-4.2}\\
&\times(\|\pa(\psi_{R} Y^{\leq 2}u^{(1)})\|_{\dot{H}^{s-1}}+\|\pa Y^{\leq 2}u^{(1)}\|_{L_x^2}) dt,\nonumber
\end{align}
for any $T\in (0,\infty)$ and
 some $C_2>\|(Y^{\leq 2}(\mu,\mu^{j}), rY^{\leq2}\mu_{0})\|_{L_{t}^{1}(L_{x}^{\infty}\cap \dot{W}^{1,n})}$.
 As $$\|\pa (\psi_{R} Y^{\leq 2}u^{(1)})(T)\|_{\dot{H}^{s-1}}
+\|\pa Y^{\leq 2}u^{(1)}(T)\|_{L_x^2}
\leq  \|u^{(1)}\|_{X_{2}}\ ,$$
by Gronwall's inequality and \eqref{ea66}, we obtain
$$\|\pa(\psi_{R} Y^{\leq 2}u^{(1)})\|_{L^\infty_t \dot{H}^{s-1}}+\|\pa Y^{\leq 2}u^{(1)}\|_{L^\infty_t L_x^2}
\leq C_{1}\ep e^{C_{2}\|(Y^{\leq 2}(\mu,\mu^{j}), rY^{\leq 2}\mu_{0})\|_{L_{t}^{1}(L_{x}^{\infty}\cap \dot{W}^{1,n})}}
\le C_1e^{C_2^2} \ep.
$$
Thus by \eqref{eqn-4.2}, we get
$$\|u^{(1)}\|_{X_{2}} \leq
C_1\ep+ C_2
\|(Y^{\leq 2}(\mu,\mu^{j}), rY^{\leq 2}\mu_{0})\|_{L_{t}^{1}(L_{x}^{\infty}\cap \dot{W}^{1,n})}C_1e^{C_2^2} \ep
\le C_3\ep\ ,$$
where we set $C_3=C_1
+C_1C_2^2e^{C_2^2}$.

If we have $\|u^{(k)}\|_{X_{2}} \leq  2 C_3\ep$ for some $k \geq 0$. Then we rewrite (\ref{4.1}) as
\beeq
\begin{cases}
\label{4.5}
\pa_{tt}u^{(k+1)} - \Delta_{\gm}u^{(k+1)} = F = F_{p}(u^{(k)}) -\mu\pa_{t}u^{(k+1)}- \mu^{j}\pa_{j}u^{(k+1)} - \mu_{0}u^{(k+1)}\\
u^{(k+1)}(0,x)=u_{0},  \pt{u^{(k+1)}}(0,x)=u_{1}.
\end{cases}
\eneq
Applying (\ref{3.4}) to (\ref{4.5}) with $F_{1} = F_{p}(u^{(k)})$, $F_{2} = -\mu\pa_{t}u^{(k+1)}- \mu^{j}\pa_{j}u^{(k+1)} - \mu_{0}u^{(k+1)}$, we get
\begin{align}
\label{4.6}
\|u^{(k+1)}\|_{X_{2}} &\leq C_{1}\ep +C_{0}\|\psi_{R}Y^{\le 2}F_2\|_{L_{t}^{1}\dot{H}^{s-1}}+C_{0}\|Y^{\le 2}F_2\|_{L_{t}^{1}L_x^{2}}\\\nonumber
&+C_{0}\|r^{-\al p}\psi_{R}^{p}Y^{\le 2}F_{p}(u^{(k)})\|_{L^{1}L^{1}L^{2}}+ C_{0}\|Y^{\le 2}F_{p}(u^{(k)})\|_{L^{1}L^{2}}.
\nonumber
\end{align}
For the norms of $F_{1}=F_{p}(u^{(k)})$, by \cite[(5.4)]{MW17} we have
$$ \|r^{-\al p}\psi_{R}^{p}Y^{\le 2}F_{p}(u^{(k)})\|_{L^{1}L^{1}L^{2}}+ \|Y^{\le 2}F_{p}(u^{(k)})\|_{L^{1}L^{2}} \les~~ \|u^{(k)}\|_{X_{2}}^{p} 
\les \ep^p\ .$$
For the norms of $F_{2} = -\mu\pa_{t}u^{(k+1)}- \mu^{j}\pa_{j}u^{(k+1)} - \mu_{0}u^{(k+1)}$, by the similar argument above we obtain
\begin{align*}
&\|\psi_{R}Y^{\le 2}F_{2}\|_{\dot{H}^{s-1}}+\|Y^{\leq 2}F_{2}\|_{ L_x^{2}}\\
\les &\|(Y^{\leq 2}(\mu,\mu^{j}), rY^{\leq 2}\mu_{0})\|_{L_{x}^{\infty}\cap \dot{W}^{1,n}}(\|\pa(\psi_{R} Y^{\leq 2}u^{(k+1)})\|_{\dot{H}^{s-1}}+\|\pa Y^{\leq 2}u^{(k+1)}\|_{L_x^2}).
\end{align*}
Hence by (\ref{4.6}) we get for any $T\in (0,\infty)$,
\begin{align}
\nonumber
&\|\pa (\psi_{R} Y^{\leq 2}u^{(k+1)})(T)\|_{\dot{H}^{s-1}}
+\|\pa Y^{\leq 2}u^{(k+1)}(T)\|_{L_x^2}\\
\leq &\|u^{(k+1)}\|_{X_{2}}
\leq  C_{1}\ep+ \CO(\ep^{p})+C_2\int_{0}^{T}\|(Y^{\leq 2}(\mu,\mu^{j}),rY^{\leq 2}\mu_{0})\|_{L_{x}^{\infty}\cap \dot{W}^{1,n}} \label{liu1}\\
&\times
(\|\pa (\psi_{R}Y^{\leq 2}u^{(k+1)})\|_{\dot{H}^{s-1}}+\|\pa Y^{\leq 2}u^{(k+1)}\|_{L_x^2}) dt\ \nonumber.
\end{align}
Applying Gronwall's inequality we have
\begin{align*}
&\|\pa (\psi_{R} Y^{\leq 2}u^{(k+1)})\|_{L^\infty_t \dot{H}^{s-1}}
+\|\pa Y^{\leq 2}u^{(k+1)}\|_{L^\infty_t L_x^2}\\
\leq& (C_{1}\ep+\CO(\ep^{p})) e^{C_{2}\|
(Y^{\leq 2}(\mu,\mu^{j}), rY^{\leq2}\mu_{0})
\|_{L_{t}^{1}(L_{x}^{\infty}\cap \dot{W}^{1,n})}}\\
\leq& (C_{1}\ep+\CO(\ep^{p})) e^{C_{2}^2}.
\end{align*}
Back to \eqref{liu1}, we obtain
\begin{align*}
\|u^{(k+1)}\|_{X_{2}} &\leq ~~C_{1}\ep+\CO(\ep^{p})+C_{2}\|\left(Y^{\leq 2}(\mu,\mu^{j}), rY^{\leq2}\mu_{0}\right)\|_{L_{t}^{1}(L_{x}^{\infty}\cap \dot{W}^{1,n})}\\
 &\times
(\|\pa(\psi_{R} Y^{\le 2}u^{(k+1)})\|_{L^{\infty}_{t}\dot{H}^{s-1}} +\|\pa Y^{\leq 2}u^{(k+1)}\|_{L_t^{\infty}L_{x}^{2}})\\
&\leq C_{1}\ep+\CO(\ep^{p})
+C_2^2(C_{1}\ep+\CO(\ep^{p})) e^{C_{2}^2}
\leq 2C_3\ep,
\end{align*}
for $0< \ep <\ep_{1}\ll 1$.
Hence $\|u^{(k)}\|_{X_{2}} \leq~~2C_{3}\ep$ for any $k \geq 0$, by induction.\par
\emph{\textbf{Convergence of the sequence}}: by (\ref{4.1}) we have
\begin{equation}
\label{4.8}
\begin{cases}
\pa_{tt}(u^{(k+1)}-u^{(k)}) - \Delta_{\gm}(u^{(k+1)}-u^{(k)}) = G,\\
(u^{(k+1)}-u^{(k)})(0,x)=0,  \pt(u^{(k+1)}-u^{(k)})(0,x)=0,
\end{cases}
\end{equation}
where $$G= F_{p}(u^{(k)})-F_{p}(u^{(k-1)})-\mu\pa_{t}(u^{(k+1)}-u^{(k)})-\mu^{j}\pa_{j}(u^{(k+1)}-u^{(k)}) - \mu_{0}(u^{(k+1)}-u^{(k)}).$$
By applying (\ref{3.4}) to (\ref{4.8}) with $G_{1} = F_{p}(u^{(k)})-F_{p}(u^{(k-1)})$ and $G_{2} = -\mu\pa_{t}(u^{(k+1)}-u^{(k)})-\mu^{j}\pa_{j}(u^{(k+1)}-u^{(k)}) - \mu_{0}(u^{(k+1)}-u^{(k)})$, then we have
\begin{align}
\label{4.9}
\|u^{(k+1)}-u^{(k)}\|_{X_{0}} \leq &C_{0}\|r^{-\al p}\psi_{R}^{p}G_{1}\|_{L^{1}L^{1}L^{2}} + C_{0}\|G_{1}\|_{L_t^{1}L_x^{2}}\\
&+C_{0}\|\psi_{R} G_{2}\|_{L_{t}^{1}\dot{H}^{s-1}}+C_{0}\|G_{2}\|_{L_t^{1}L^{2}_{x}}.\nonumber
\end{align}
By \cite[(5.5)]{MW17}, we have
\begin{eqnarray*}
\|r^{-\al p}\psi_{R}^{p}G_{1}\|_{L^{1}L^{1}L^{2}} + \|G_{1}\|_{L_t^{1}L_x^{2}}& \les & \|(u^{(k)},u^{(k-1)})\|_{X_{2}}^{p-1}\|u^{(k)}-u^{(k-1)}\|_{X_{0}} \\
 & \les & \ep ^{p-1}\|u^{(k)}-u^{(k-1)}\|_{X_{0}}.
\end{eqnarray*}
For $G_{2}$ part, we have
\begin{align*}
&\|\psi_{R} G_{2}\|_{ \dot{H}^{s-1}}+\|G_{2}\|_{ L^{2}_{x}}\\
\les& \|(\mu,\mu^{j}, r\mu_{0})\|_{L_{x}^{\infty}\cap \dot{W}^{1,n}}(\|\pa(\psi_{R}(u^{(k+1)}-u^{(k)}))\|_{\dot{H}^{s-1}}+\|\pa(u^{(k+1)}-u^{(k)})\|_{L^{2}_x}).
\end{align*}
Thus by (\ref{4.9}), we have for any $T\in(0,\infty)$,
\begin{align}
&\|\pa\psi_{R}(u^{(k+1)}-u^{(k)})(T)\|_{\dot{H}^{s-1}}+\|\pa(u^{(k+1)}-u^{(k)})(T)\|_{L^{2}_{x}}\nonumber\\
 \leq &\|u^{(k+1)}-u^{(k)}\|_{X_{0}}\label{eqn-4.3}\\
\leq &C_{4} \ep ^{p-1}\|u^{(k)}-u^{(k-1)}\|_{X_{0}} +C_{4}\int_{0}^{T}\|(\mu,\mu^{j}, r\mu_{0})\|_{L_{x}^{\infty}\cap \dot{W}^{1,n}}\nonumber\\
 & \times(\|\pa\psi_{R}(u^{(k+1)}-u^{(k)})\|_{\dot{H}^{s-1}}+\|\pa(u^{(k+1)}-u^{(k)})\|_{L^{2}_{x}})dt,\nonumber
\end{align}
for some $C_4>C_2>0$.
By Gronwall's inequality, for any $T < \infty$, we have
$$\|\pa\psi_{R}(u^{(k+1)}-u^{(k)})(T)\|_{\dot{H}^{s-1}}+\|\pa(u^{(k+1)}-u^{(k)})(T)\|_{L^{2}_{x}} \leq C_{4}\ep^{p-1}
e^{C_4^2}\|u^{(k)}-u^{(k-1)}\|_{X_{0}} .$$
Hence by \eqref{eqn-4.3}
\begin{eqnarray*}
\|u^{(k+1)}-u^{(k)}\|_{X_{0}} & \le & C_{4} \ep ^{p-1}\|u^{(k)}-u^{(k-1)}\|_{X_{0}}+
C_{4}^3 \ep^{p-1}
e^{C_4^2}\|u^{(k)}-u^{(k-1)}\|_{X_{0}}
 \\
 & \le & \frac{1}{2}\|u^{(k)}-u^{(k-1)}\|_{X_{0}},
\end{eqnarray*}
for $0 <\ep< \ep_{2}$ with $\ep_2\ll 1$,
which yields the sequence converges in $X_{0}$. When $\ep_0=\min(\ep_1,\ep_2)$ and $\ep\in (0,\ep_0)$, the limit $u \in X_{2}$ with $\|u\|_{X_{2}} \les~~\ep$ is the solution we are looking for.
\end{prf}
\section{Almost Global Existence}
In this section, we prove the almost global existence for the four-dimensional critical problem,
Theorem \ref{main.1}.
We set, for $T\in (0,\infty)$,
\begin{align*}
\|u\|_{\tilde{X}_T^{m}} := &(\ln T)^{-1/4}\|\psi_{R}\<r\>^{-\frac{1}{4}} Y^{\leq m}u\|_{L_T^{4}L^{2}}+\|Y^{\leq m}u\|_{LE_{T}}\\
&+\|\langle r\rangle^{-\frac{3}{4}}Y^{\leq m}u\|_{L_{T}^{2}L^{2}}
+\|\pa (\psi_{R}Y^{\leq m}u)\|_{L_{T}^{\infty}\dot{H}^{-1}}.
\end{align*}

\begin{theorem}
\label{almost}
Let $n =4$, and assume that \eqref{ea6}, \eqref{ea8} hold. Consider \eqref{1.1} with \eqref{6666} and $F_{p}= u^{2}$. Then there exist $\ep_{0}>0$ sufficiently small, $R$ sufficiently large and $c > 0$, so that if $0 < \ep < \ep_{0}$ and
\beeq
\label{eqn-5.1}
\|\nabla^{\le 1}Y^{\le 3}u_{0}\|_{L^{2}_{x}}+ \|Y^{\le 3}u_{1}\|_{L^{2}_{x}}+\|\psi_{R}Y^{\le 3}u_{1}\|_{\dot{H}^{-1}} \leq \ep
\eneq
then there is a unique solution $u \in [0, T_\ep] \times \R^{4}$ with
$\|u\|_{\tilde{X}^{3}_{T_\ep}}\les \ep$,
where
$T_\ep=\exp(c\ep^{-2})$.
\end{theorem}
\begin{prf}
Basically, the proof follows the similar way in Theorem \ref{global}.
For convenience of statement, we introduce
\begin{align*}
\|F\|_{\tilde{N}_T^{m}} = \|\psi_{R}^{2} r^{-1/2}Y^{\leq m}F_{1}\|_{L_{T}^2L^1H^{-1/2+}_{\omega}}
+\|Y^{\leq m}(F_{1},F_{2})\|_{L^{1}_{T}L^{2}}+\|\psi_{R}Y^{\leq m}F_{2}\|_{L_{T}^{1}\dot{H}^{-1}}.
\end{align*}
Then by Lemma \ref{thm-wStri-homo} with $n = 4$, for linear equation (\ref{3.1}), there exists a constant $C_{5} > 0$ such that for any $0 \leq m \leq 3$, we have
\beeq
\label{5.2}
\|u\|_{\tilde{X}_T^{m}} \leq C_{5}\|\pa^{\le 1} Y^{\leq m}u(0,\cdot)\|_{L^{2}}+C_{5}\|\psi_{R}Y^{\leq m}\pa_{t}u(0)\|_{\dot{H}^{-1}}+ C_{5}\|F\|_{\tilde{N}^{m}_T}.
\eneq
We set $u^{(0)} =0$ and recursively define $u^{(k+1)}$ be the solution to the linear equation
\begin{equation}
\label{5.1}
\begin{cases}
\pa_{tt}u^{(k+1)} - \Delta_{\gm}u^{(k+1)} +
 \mu \pt u^{(k+1)} +
 \mu^{j}\pa_{j}u^{(k+1)} + \mu_0u^{(k+1)}  = (u^{(k)})^{2},\\
u^{(k+1)}(0,x)=u_{0},  \pt{u^{(k+1)}}(0,x)=u_{1}.
\end{cases}
\end{equation}
To complete the proof, we need to show that $u^{(k)}$ is well defined, bounded in $\tilde{X}_{T_\ep}^{3}$ and convergent in $\tilde{X}_{T_\ep}^{0}$.

\emph{\textbf{Well defined}}: It is easy to see $u^{(1)} \in C([0,\infty);H^{4})\cap C^{1}([0,\infty);H^{3})$. If we have $u^{(k)} \in C([0,\infty);H^{4})\cap C^{1}([0, \infty);H^{3})$ for some $k\geq 1$.
Since $$|\nabla^{\leq 3}(u^{(k)})^{2}| \lesssim |\nabla^{\leq1}u^{(k)}||\nabla^{\leq2}u^{(k)}| + |u^{(k)}\nabla^{\leq3}u^{(k)}|,$$
by H\"{o}lder's inequality and Sobolev embedding, we have, for any $t\in [0,T]$ with fixed $T<\infty$,
\begin{align*}
\|\nabla^{\leq 3}(u^{(k)})^{2}\|_{L^{2}} &\les \|\nabla^{\leq1}u^{(k)}\|_{L^{4}}\|\nabla^{\leq2}u^{(k)}\|_{L^{4}}+ \|u^{(k)}\|_{L^{\infty}}\|\nabla^{\leq3}u^{(k)}\|_{L^{2}}\\
&\les\|\nabla^{\leq1}u^{(k)}\|_{\dot{H}^{1}}\|\nabla^{\leq2}u^{(k)}\|_{\dot{H}^{1}}+\|u^{(k)}\|_{H^{3}}
\|\nabla^{\leq3}u^{(k)}\|_{L^{2}} < \infty.
\end{align*}
Thus $(u^{(k)})^{2} \in L^{1}_{t}([0,T];H^{3})$ for any $T < \infty$. By standard existence theorem of linear wave equations, we have $u^{(k+1)} \in C([0,T];H^{4})\cap C^{1}([0,T];H^{3})$  for any $T < \infty$ and so
$u^{(k+1)} \in C([0,\infty);H^{4})\cap C^{1}([0, \infty);H^{3})$. Hence the iteration sequence is well defined.

In the following, we give the proof of boundedness and omit the similar proof of convergence.\\
\emph{\textbf{Boundedness in $\tilde{X}_T^{3}$}}: As usual, we prove the boundedness by induction. Assuming $\|u^{(k)}\|_{\tilde{X}_T^{3}} \leq 2C_{6}\ep $, for some $k \geq 0$ and $C_{6}$ to be determined later. Then we rewrite (\ref{5.1}) as
\begin{equation}
\label{5.3}
\begin{cases}
\pa_{tt}u^{(k+1)} - \Delta_{\gm}u^{(k+1)} = F =(u^{(k)})^{2}- \mu \pt u^{(k+1)} -
 \mu^{j}\pa_{j}u^{(k+1)} - \mu_0u^{(k+1)},\\
u^{(k+1)}(0,x)=u_{0},  \pt{u^{(k+1)}}(0,x)=u_{1}.
\end{cases}
\end{equation}
By applying (\ref{5.2}) to (\ref{5.3}) with $F_{1} = (u^{(k)})^{2}$ and $F_{2}= - \mu \pt u^{(k+1)} -
 \mu^{j}\pa_{j}u^{(k+1)} - \mu_0u^{(k+1)}$, we obtain
\beeq
\label{5.5}
\|u^{(k+1)}\|_{\tilde{X}_T^{3}} \leq C_{5}\ep+ C_{5}\|F\|_{\tilde{N}^{3}_T}.
\eneq
For the norms of $F_{1}$ part, by \cite[(5.15)]{W17},  we have
$$\|\psi_{R}^{2} r^{-1/2}Y^{\leq 3}(u^{(k)})^{2}\|_{L_{T}^2L^1H^{-1/2+}_{\omega}}
+\|Y^{\leq 3}(u^{(k)})^{2}\|_{L^{1}_{T}L^{2}} \les (\ln T)^{1/2}\|u^{(k)}\|^{2}_{\tilde{X}_T^{3}}\ ,$$
for any $2 \leq T <\infty$.
For the norms of $F_{2}$ part, by the similar argument of \eqref{eq-4.1}-\eqref{liu7}, we have
\begin{align*}
&\|\psi_{R}Y^{\le 3}F_{2}\|_{\dot{H}^{-1}}+\|Y^{\leq 3}F_{2}\|_{ L_x^{2}}\\
\les &\|(Y^{\leq 3}(\mu,\mu^{j}), rY^{\leq 3}\mu_{0})\|_{L_{x}^{\infty}\cap \dot{W}^{1,4}}(\|\pa(\psi_{R} Y^{\leq 3}u^{(k+1)})\|_{\dot{H}^{-1}}+\|\pa Y^{\leq 3}u^{(k+1)}\|_{L_x^2}).
\end{align*}
Thus by (\ref{5.5}) we get for any $t\in [0, T]$ and $2 \leq T <\infty$,
\begin{align}
\nonumber
&\|\pa (\psi_{R}Y^{\leq 3}u^{(k+1)})(t)\|_{\dot{H}^{-1}}
+\|\pa Y^{\leq 3}u^{(k+1)}(t)\|_{L^{2}}\\
\leq & \|u^{(k+1)}\|_{\tilde{X}_T^{3}}
\label{liu5}\\
\leq&  C_{5}\ep + C_{7}(\ln T)^{1/2}\|u^{(k)}\|^{2}_{\tilde{X}_T^{3}}
+C_{7}\int_{0}^{t}\|\left(Y^{\leq 3}(\mu,\mu^{j}), rY^{\leq 3}\mu_{0}\right)\|_{L_{x}^{\infty}\cap \dot{W}^{1,4}} \nonumber\\
&\times
(\|\pa(\psi_{R} Y^{\leq 3}u^{(k+1)})\|_{\dot{H}^{-1}}+\|\pa Y^{\leq 3}u^{(k+1)}\|_{L_x^2})d\tau\ \nonumber,
\end{align}
for some $C_{7} > \|(Y^{\leq 3}(\mu,\mu^{j}), rY^{\leq 3}\mu_{0})\|_{L_{t}^{1}(L_{x}^{\infty}\cap \dot{W}^{1,4})}$ and $C_{7}$ is independent of $C_{6}$ and $T$.
By Gronwall's inequality, for any $t\in [0, T]$, we obtain
$$\|\pa (\psi_{R}Y^{\leq 3}u^{(k+1)})(t)\|_{\dot{H}^{-1}}
+\|\pa Y^{\leq 3}u^{(k+1)}(t)\|_{L^{2}} \leq (C_{5}\ep + C_{7}(\ln T)^{1/2}\|u^{(k)}\|^{2}_{\tilde{X}_T^{3}})e^{C_{7}^{2}},$$
and by (\ref{liu5}) again,
\begin{align*}
\|u^{(k+1)}\|_{\tilde{X}_T^{3}} \leq C_{5}\ep + C_{7}(\ln T)^{1/2}\|u^{(k)}\|^{2}_{\tilde{X}_T^{3}} + C^{2}_{7}(C_{5}\ep + C_{7}(\ln T)^{1/2}\|u^{(k)}\|^{2}_{\tilde{X}_T^{3}})e^{C_{7}^{2}}.
\end{align*}
If we take $C_{6} = C_{5}+ C_{5}C^{2}_{7}e^{C_{7}^{2}}$ and $T = T_{\ep} = e^{c\ep^{-2}}$ for some $c > 0$ such that
$$4c^{1/2}C_{6}C_{7}(1+ C_{7}^{2}e^{C_{7}^{2}}) \leq 1,$$
then
\begin{align*}
\|u^{(k+1)}\|_{\tilde{X}^{3}_{T_{\ep}}} \leq C_{6}\ep +4c^{1/2}C^{2}_{6}C_{7}(1+ C_{7}^{2}e^{C_{7}^{2}})\ep \leq 2C_{6}\ep.
\end{align*}
Hence $\|u^{(k)}\|_{\tilde{X}^{3}_{T_{\ep}}} \leq 2C_{6}\ep$ for any $k \geq 0$, by induction.
\end{prf}

\end{document}